\documentclass[11pt,a4paper]{article}
  
\usepackage{latexsym}
\usepackage{graphicx}
\usepackage{helvet}
\usepackage{psfrag}
\usepackage{bbm}
\usepackage{bm}
\usepackage{setspace}
\usepackage{comment}

\usepackage{natbib}
\usepackage{amsmath}  
\usepackage{mathabx}  
\usepackage{multirow}
\usepackage{url}

\renewcommand{\vec}{\bm}

\newcommand{\floor}[1]{\left\lfloor #1\right\rfloor}
\newcommand{\ceil}[1]{\left\lceil #1\right\rceil}

\def\Chi{%
    \mbox{
    {\kern-.20em\setbox0=\hbox{X}%
    \vbox to 1.0\ht0{\hbox{$\chi$}\vss}}%
    \kern-.00em} }
\def\mathlatex{%
    \mbox{
    L\kern-.36em {\setbox0=\hbox{T}%
    \vbox to \ht0{\hbox{\the\scriptfont0 A}\vss}}%
    \kern-.15em \TeX} }

\newcommand{\qed}{\hbox{}\hfill $\Box$ }

\newenvironment{blinded}{}{}

\newenvironment{keywords}{\noindent\textbf{Keywords:}}{}
\renewcommand{\and}{;}

\excludecomment{twocolversion}

\newcommand{\papertitle}{Constructing wavelets by welding segments of smooth
functions}

\title{\papertitle}
\begin{blinded}
\author{Maarten Jansen\\
Universit\'e libre de Bruxelles, departments of Compter Science and Mathematics}
\end{blinded}

\begin{document}

\maketitle

\begin{abstract}
The construction of B-spline wavelet bases on nonequispaced knots is
extended to wavelets that are piecewise segments from any combination of smooth
functions. 
The extended wavelet family thus provides multiresolution basis functions
with support as compact as possible and belonging to a user controlled
smoothness class.
The construction proceeds in two phases. In the first fase, a set of smooth
functions is used in the welding of compact supported, piecewise
smooth basis functions. These piecewise smooth basis functions are refinable,
meaning that they can be written as linear combinations of similar basis
functions constructed on a fined grid of knots. The expression of the linear
combination between the bases at two scales is known as a refinement
or two-scale equation.
In the second phase, the refinabability enables the construction of a wavelet
transform. To this end, the refinement equation of the piecewise smooth scaling
functions is factored into a lifting scheme, to which the desired properties 
of the subsequent wavelet basis can then be added.
Next to the details of the construction, the paper discusses the conditions for
it to fit into the classical framework of multiresolution analyses.

\end{abstract}
\newpage

\section{Introduction}

The construction of a wavelet basis and corresponding transform may proceed in
either of two directions. The transform first approach starts by constructing a
perfect reconstruction (i.e., invertible) filterbank which is then applied in
an iterative way to define a multiresolution transform. The corresponding basis
then follows from subdivision or infinite refinement, which is essentially an
inverse transform up to infinitely fine resolution level. The analysis of the
smoothness properties of the resulting basis functions is, however, notoriously
difficult, especially when the refinement takes place on nonequispaced knots
\citep{daubechies99:regir,daubechies01:commutation}.
On the positive side, the lifting scheme provides a straightforward tool for
the construction of invertible filterbanks with desired properties
\citep{sweldens96:liftingbior}.
These properties may include the exact reproduction of polynomials, meaning
that the space spanned by the basis includes all polynomials up to a given
degree.
This property is termed dual vanishing moments in the context of bi-orthogonal
wavelet transforms. Another feature of interest, especially in statistical
applications, is the control of variance propagation through carful design of
the lifting scheme
On nonequispaced knots, the Deslauriers-Dubuc 
\citep{deslauriers89:symmetric,donoho99:deslauriers} refinement adopts
polynomial interpolation in the design of a lifting scheme. The polynomial
interpolation ensures the above mentioned polynomial reproduction after
subdivision.

The basis first approach starts by looking for a refinable basis. A basis
defined on a given set of knots is called refinable if it can be written as a
linear combination of similar basis functions defined on set of the same and
additional knots, thus creating a sequence of nested sets of knots.
A typical example of refinable bases are B-splines \citep{deBoor01:splines}.
Splines of order $\widetilde{\nu}$ are piecewise polynomials of degree
$\widetilde{\nu}-1$ on a given set of knots. That means that each segment of
the spline function between two adjacent knots coincides with a polynomial.
Moreover, the transitions from one segment to another in the knots have
$(\widetilde{\nu}-2)$ continuous derivatives. Finally, B-splines are defined to
have the smallest possible support, being nonzero on $\widetilde{\nu}$
adjacent intervals bounded by the knots. B-splines constitute a basis for all
splines, including the plain polynomials of degree $\widetilde{\nu}-1$.
As a result, a B-spline defined by a coarse set of knots is a linear
combination of the B-splines defined on any fine set of knots that includes
the given coarse set.
When this refinement relation is factored into lifting steps
\citep{daubechies98:factorlifting}, the procedure yields a primitive wavelet
transform. A final lifting step is then sufficient to construct any possible
wavelet transform for this refinement \citep{jansen22:waveletbook}.
The basis first approach applied to B-splines on nonequispaced, nested knots
extends the dyadic Cohen-Daubechies-Feauveau spline wavelets
\citep{coh-dau-fea92:biorthogonal}.

This paper further extends the nonequispaced B-spline wavelets by replacing the
polynomials as smooth building blocks by any other collection of
$\widetilde{\nu}$ smooth functions. 
In Section \ref{sec:brokenbasis} a broken basis on compact intervals is
constructed starting from non-local smooth functions.
In Section \ref{sec:BBrefinement}, the refinement or two-scale equation 
for a broken basis is developed. This expression is the key to a multiscale
transformation.
Section \ref{sec:BBwavelets} then adds the detail or wavelet basis, so that in
Section \ref{sec:BBWT} the basis transform, i.e., the wavelet transform can be
completed.
Section \ref{sec:discussion} investigates the properties of the proposed
wavelet basis, finding out that on equidistant sets of knots, it generally
lacks the classical property of being dilations and translations of a single
father and mother function. It is explained that this may have an important
impact, even on nonequidistant knots, where a father and mother functions are
not feasible anyway.
The paper is concluded in the short Section \ref{sec:conclusion}.

\section{The construction of a broken basis }
\label{sec:brokenbasis}

The objective of this section is to present a general technique for the
construction of basis functions with compact support consisting of segments
from smooth, long-ranging functions, with smooth transitions between the
segments in the knots. This way we build basis functions that are local in
space (or time), exactly as B-splines do starting from power functions. In the
next section, a multiscale, wavelet basis will add locality in scale to the
spatial locality of this section.

\subsection{Setting up the equations in the knots}

Let $\Omega(x) = \left[\omega_0(x) \ldots \omega_{\widetilde{\nu}-1}(x)\right]$
be a row vector of linearly independent $(\widetilde{\nu}-1)$ continuously
differentiable functions $\omega_q(x)$, defined on the interval
$[0,1]$. In the case of B-splines the functions in $\Omega(x)$ would be the
power functions, i.e., $\omega_q(x) = x^q$. Furthermore, define the nested sets 
of knots $X_j = \{x_{j,0},x_{j,1},\ldots,x_{j,n_j-1}\}$ with
$X_j \subset X_{j+1}$ and $n_j$ the number of knots at resolution level $j$.
Although not strictly necessary for the construction, the left most knot is
taken to be $x_0=0$, while the right most knot is taken to be $x_{n_j-1}=1$.
Then the broken basis at level $j$ is defined as the row vector $\Phi_j(x)$ of
$n_j+\widetilde{\nu}-2$ basis functions
\[
\Phi_j(x) =
\left[
\varphi_{j,0}(x)\,\varphi_{j,1}(x)\,\ldots\,\varphi_{j,n_j+\widetilde{\nu}-3}(x)
\right],
\]
where each basis function consists of a number of segments,
\begin{equation}
\varphi_{j,k}(x) =
\sum_{i\in S_{j,k}} \chi_{j,i}(x) \Omega(x) \vec{a}_{j,k,i}.
\label{eq:defBB}
\end{equation}
In this expression, $\chi_{j,i}(x)$ is the indicator function on the 
subinterval $[x_{j,i},x_{j,i+1}]$,
and the index set $S_{j,k}$ is given by
$S_{j,k} = \{0,1,\ldots,n_j-2\} \cap \{k-\widetilde{\nu}+1,\ldots,k\}$.
The support of $\varphi_{j,k}(x)$ is then given by $[x_{j,l(k)},x_{j,r(k)}]$,
where $l(k) = \min S_{j,k}$ and $r(k) = \max S_{j,k}+1$.
It should be noted that all indexing of sets, matrices and vectors throughout
this text starts with zero. Hence, for instance, the first
element of a matrix $\mathbf{A}$ will be $A_{0,0}$.

The coefficient vectors $\vec{a}_{j,k,i}$ follow from the continuity and
boundary conditions. For $r=0,1,\ldots,\widetilde{\nu}-2$ it is imposed that in
all interior knots $i \in \{1,2,\ldots,n_j-2\}$,
\begin{equation}
\lim_{x \to x_{j,i}^-} {d^r \over dx^r} \Phi_j(x) =
\lim_{x \to x_{j,i}^+} {d^r \over dx^r} \Phi_j(x).
\label{eq:BBknotcontinuity}
\end{equation}
In the left boundary, we impose that for $r \neq k$,
\begin{equation}
\lim_{x \to x_{j,0}^+} {d^r \over dx^r} \varphi_{j,k}(x) = 0.
\label{eq:BBleftboundarycond}
\end{equation}
In the right boundary, a similar condition holds for $r \neq
n_j+\widetilde{\nu}-3-k$,
\begin{equation}
\lim_{x \to x_{j,n_j-1}^-} {d^r \over dx^r} \varphi_{j,k}(x) = 0.
\label{eq:BBrightboundarycond}
\end{equation}
Finally, the basis functions are normalised by imposing
\begin{equation}
\sum_{k=0}^{n_j+\widetilde{\nu}-3} \varphi_{j,k}(x) = \omega_0(x).
\label{eq:BBpartitionofunity}
\end{equation}
In other words, knowing that all functions in $\Omega(x)$ are linear
combinations of $\Phi_j(x)$, i.e., $\Omega(x) = \Phi_j(x) \mathbf{C}_j$, we
impose that the first column of $\mathbf{C}_j$ are all ones. If $\omega_0(x) =
1$ (which is a typical choice), then this condition is known as the
\emph{partition of unity}.

\subsection{Solving the equations}

Conditions (\ref{eq:BBknotcontinuity}) to (\ref{eq:BBpartitionofunity}) provide
linear equations in the unknown coefficient vectors $\vec{a}_{j,k,i}$. 
As each vector contains $\widehat{\nu}$ coefficients, the total number of
unknowns follows the number of vectors, which is given for each basis function
by the number of segments. Of the $n_j+\widetilde{\nu}-2$, the two boundary
functions each have one segment, the next adjacent have two segments, and so
on, up to a maximum of $\widetilde{\nu}$ segments. The number of basis
functions with $\widetilde{\nu}$ segments is given by
$n_j+\widetilde{\nu}-2-2(\widetilde{\nu}-1) = n_j-\widetilde{\nu}$.
All together, we find 
\(
(n_j-\widetilde{\nu}-2) \cdot \widetilde{\nu} +
(1+2+\ldots+\widetilde{\nu}-1) \cdot 2 =
(n_j-\widetilde{\nu})\widetilde{\nu}+\widetilde{\nu}(\widetilde{\nu}-1)
=
(n_j-1)\widetilde{\nu}
\)
vectors, hence
\(
(n_j-1)\widetilde{\nu}^2
\)
unknown coefficients.

As for the number of equations, the continuity expressions in
(\ref{eq:BBknotcontinuity}) lead to $\widetilde{\nu}-1$ equations in each knot
and for each basis function with nonzero value or derivatives in that knot.
The fully interior basis functions, with support in $[x_{j,1},x_{j,n_j-2}]$
that is, are completely defined by these continuity expressions, up to the
normalisation in (\ref{eq:BBpartitionofunity}).
Each of these functions has $\widetilde{\nu}+1$ knots with
$\widetilde{\nu}-1$ equations. The number of these functions is
given by $n_j-2-\widetilde{\nu}$, leading to
$(n_j-2-\widetilde{\nu})(\widetilde{\nu}^2-1)$ equations. For $k =
0,1,\ldots,\widetilde{\nu}-1$, the functions $\varphi_{j,k}(x)$ and
$\varphi_{j,n_j+\widetilde{\nu}-3-k}(x)$ are on the boundary.
They have $k+1$ interior knots, each with $\widetilde{\nu}-1$ equations and one
boundary point with the same number of equations. The boundaries thus add a
total of $(\widetilde{\nu}^2-1)\widetilde{\nu}$ equations. Finally, the
normalisation in (\ref{eq:BBpartitionofunity}) is imposed in $n_j-1$ intervals
$[x_{j,i},x_{j,i+1}]$, leading to a total of
\(
(n_j-2-\widetilde{\nu})(\widetilde{\nu}^2-1) +
(\widetilde{\nu}^2-1)\widetilde{\nu} + (n_j-1)
=
(n_j-1) \widetilde{\nu}^2
\)
linear equations, uniquely solving the unknown coefficients.

It is interesting to see that neither the values of the knots, $x_{j,k}$, nor
the function values of $\Phi_j(x)$ or $\Omega(x)$ outside the knots, have any
direct impact on the solution of the system for the break up coefficients
$\vec{a}_{j,k,i}$.
Of course, the knots $x_{j,k}$ have indirect influence, through the function
values $\Phi_j(x_{j,k})$ and the derivatives evaluated in these knots, but
given these values, the coefficients $\vec{a}_{j,k,i}$ follow in a unique way.
Conversely, the coefficients $\vec{a}_{j,k,i}$ are enough to find the entries
of the $(n_j+\widetilde{\nu}-2)\times \widetilde{\nu}$ matrix 
$\mathbf{A}_j$ in the expansion
\(
\Omega(x) = \Phi_j(x)\mathbf{A}_j.
\)
Indeed, on each subinterval $[x_{j,i},x_{j,i+1}]$, consider the coefficient
vectors $\vec{a}_{j,k,i}$ corresponding to broken basis functions
$\varphi_{j,k}(x)$ that are nonzero on that interval, i.e.,
\(
k = \{i,i+1,\ldots,i+\widetilde{\nu}-1\}.
\)
Then the submatrix $\mathbf{A}_{j;i}$ consisting of rows
$(i,i+1,\ldots,i+\widetilde{\nu}-1)$ follows from
\[
\mathbf{A}_{j;i} = \left[\begin{array}{cccc}
\vec{a}_{j;i,i} & \vec{a}_{j;i+1,i} & \ldots & \vec{a}_{j;i+\widetilde{\nu}-1,i}
\end{array}\right]^{-1}.
\]
As the submatrices $\mathbf{A}_{j;i}$ have overlapping columns, finding the
elements of $\mathbf{A}_j$ is an overdetermined problem. Taking averages of
solutions reduces the numerical errors.

\subsection{Interpolation between the knots}

As the evaluations of $\Omega(x)$ and its derivatives in the knots are
sufficient for the construction of the broken basis $\Phi_j(x)$, Hermite
polynomial interpolation (i.e., interpolation of function values and
derivatives) may be used in each subinterval $[x_{j,i},x_{j,i+1}]$ separately
to define values of $\Omega(x)$ and $\Phi_j(x)$ outside the knots, at least if
no such values are given already. 
Interpolation on each subinterval leads to the approximation of
$\Omega(x)$ and $\Phi_j(x)$ by Hermite splines of order $2\widetilde{\nu}$,
i.e., piecewise polynomials of degree $2\widetilde{\nu}-1$ defined by the
function values and $\widetilde{\nu}-1$ derivatives in the two
end points $x_{j,i}$ and $x_{j,i+1}$ of the interval, where all derivatives
except for the $(\widetilde{\nu}-1)$st are continuous.
The approximations are useful in computing the values of the moments
\[
M_{j;k,i,q} = \int_{x_{j,i}}^{x_{j,i+1}} \varphi_{j,k}(x) x^q \, dx,
\]
which is of interest in the construction of wavelets, as discussed in Section
\ref{sec:BBwavelets}.

\begin{figure}[t!]
\hspace*{-20mm}
\includegraphics[width=1.2\textwidth]{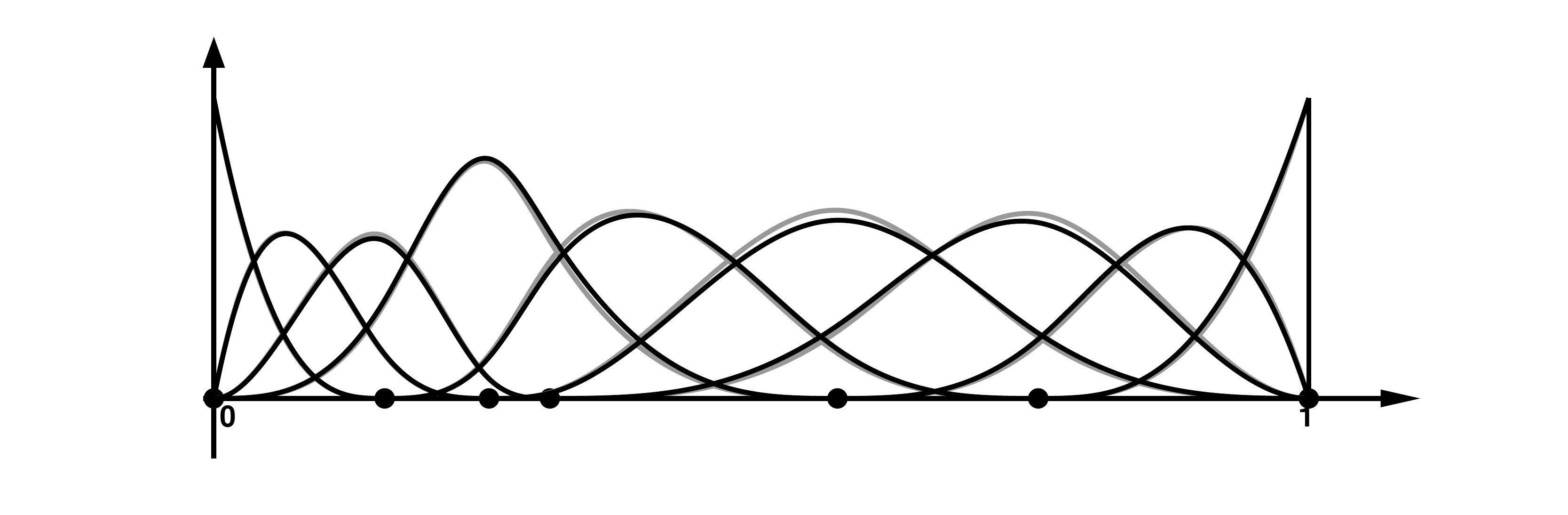}
\\
\hspace*{-20mm}
\includegraphics[width=1.2\textwidth]{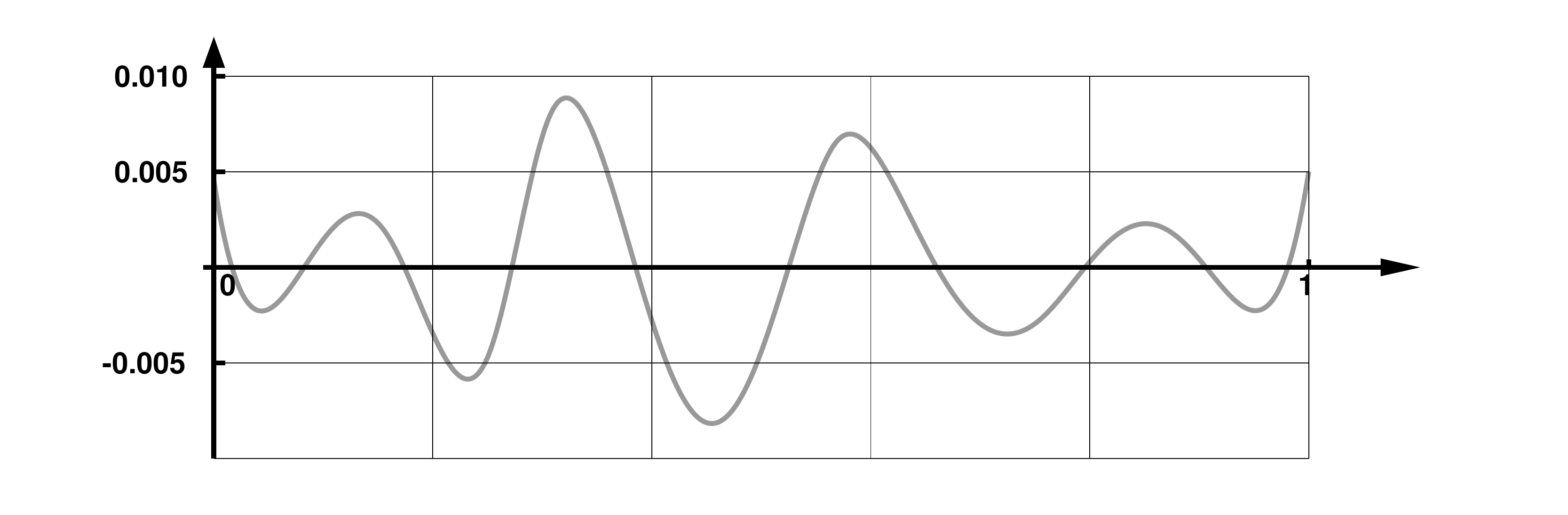}
\caption{Refinable broken bases formed from different smooth bases.
Top panel, in grey line:
cubic B-splines, formed from
$\Omega(x) = \left[\begin{array}{cccc}1 & x & x^2 & x^3 \end{array}\right]$.
In black line: basis formed with segments from
$\Omega(x) = \left[\begin{array}{cccc}1 & x & \sin(2\pi x) & \cos(2\pi x)
\end{array}\right]$.
Bottom panel, in grey line: approximation error curve of the orthogonal
projection of $y(x) = \cos(2\pi x)$ onto the cubic B-spline basis.
The error curve in the basis with sine and cosine segments is zero up to
rounding errors.
}
\label{fig:brokenbases}
\end{figure}

Figure \ref{fig:brokenbases} illustrates the construction of refinable broken
bases on the interval $[0,1]$.
It depicts two different bases. The functions in grey line are cubic B-splines
on a set of $n_j=7$ nonequispaced knots, where $x_0=0$ and $x_6 = 1$.
Cubic B-splines are constructed from
$\Omega(x) = \left[\begin{array}{cccc}1 & x & x^2 & x^3 \end{array}\right]$.
Replacing $x^2$ and $x^3$ by $\sin(2\pi x)$ and $\cos(2\pi x)$ leads to the
construction of the basis in black lines.
The two bases are quite similar. The black basis offers the advantage of a
perfect reconstruction of the $\sin(2\pi x)$ and $\cos(2\pi x)$ functions,
whereas in the grey basis the orthogonal projection of these functions induce
an pointwise approximation error of the order of magnitude of $0.01$.

\subsection{Local specification of smooth functions}

The construction of $\Phi_J(x)$ at fine resolution $j$ from a fixed smooth
$\Omega(x)$ may suffer from numerical ill conditioning. 
Indeed, let $h_J$ be the fine working scaling, then
\(
\omega_q(x_{J,k}+h_J) \approx \omega_q(x_{J,k}) + \omega_q'(x_{J,k}) h_J,
\)
leading to a near linear dependence in the building set of basis functions.
The ill condition can be dealt with by choosing local linear combinations
$\Omega_{J,k}(x) = \Omega(x)_{J,k}A_{J;k}$ for the construction
of $\varphi_{J,k}(x)$, assuming that such local basis functions can be
specified in a numerically well conditioned way
(in particular, if possible, without explicit calculation of $A_{J;k}$).
As an example, a building set of power functions with
$\omega_q(x) = x^q$ can be replaced by a local working set of locally centered
power functions $\omega_{J,k,q}(x) = (x-x_{J,k})^q$, where both constructions
would lead to a B-spline basis.

The local specification also enables the assembling of pieces from different
smooth functions on subsequent segments, thus adding an almost infinite
flexibility to the design. Applications in the numerical solution of
differential equations under intermittent conditions are obvious.

\section{Broken basis refinement}
\label{sec:BBrefinement}

By construction, the basis in (\ref{eq:defBB}) is refinable.
More precisely, suppose there exists a grid $\vec{x}_{j+1}$ of knots
$x_{j+1,k}$ with $k \in \{0,1,\ldots,n_{j+1}\}$ and
$n_{j+1}\times n_j$ matrix $\widetilde{\mathbf{J}}_j$ with elements in
$\{0,1\}$ so that $\vec{x}_j = \widetilde{\mathbf{J}}_j^T \vec{x}_{j+1}$, i.e.,
$\vec{x}_j$ is obtained by subsampling from $\vec{x}_{j+1}$.
With these nested sets of knots refinability means that there also exists a
$(n_{j+1}+\widetilde{\nu}-2)\times(n_j+\widetilde{\nu}-2)$ matrix
$\mathbf{H}_j$ so that a two-scale equation holds,
\begin{equation}
\Phi_j(x) = \Phi_{j+1}(x) \mathbf{H}_j.
\label{eq:2scale}
\end{equation}
The refinable functions in $\Phi_j(x)$ and $\Phi_{j+1}(x)$ are also referred to
as scaling functions.
The remainder of this section is devoted to finding the entries of the
refinement matrix $\mathbf{H}_j$.
In the subsequent discussion, the set
$O_{j+1} \subset \{0,1,\ldots,n_{j+1}-1\}$
will denote the indices $l$ for which $x_{j+1,l}$ is not a component of
$\vec{x}_j$.
Typically the set $O_{j+1}$ contains every second index,
$O_{j+1} = \{1,3,5,\ldots\}$.
The set will be referred to as the ``odd indices'', even if partitionings 
other than even-odd remain possible.
The complement of $O_{j+1}$ in $\{0,1,\ldots,n_{j+1}-1\}$ are the
``even indices'', again in broad sense, and written as $E_{j+1}$.
The vectors of all knots with respectively even and odd indices at level $j+1$
are denoted by $\vec{x}_{j+1,e}$ and $\vec{x}_{j+1,o}$.

The subsequent discussion concentrates on the even-odd partitioning, i.e., the
case where $O_{j+1} = \{1,3,5,\ldots\}$.
In that case, it is straightforward to prove that the columns of the matrix
$\mathbf{H}_j$ have at most $\widetilde{\nu}+1$ nonzeros, see for instance
\citep[Corollary 4.2.16, page 117]{jansen22:waveletbook}.
In general, $H_{j;\ell,k}$ can only be nonzero if the support of
$\varphi_{j+1,\ell}(x)$ is a subset of the support of $\varphi_{j,k}(x)$.
In order to develop this condition, we need to deal with the boundaries.
For this reason, we define knots with negative index to coincide with the left
boundary $x_{j,0} = x_{j+1,0}$ and knots with index above $n_j$ or $n_{j+1}$ to
coincide with the right boundaries $x_{j,n_j}$ and $x_{j+1,n_{j+1}}$.
Then, using (\ref{eq:defBB}), we find
\(
H_{j;\ell,k} \neq 0 
\Rightarrow
[x_{j+1,\ell-\widetilde{\nu}+1},x_{j+1,\ell+1}] \subset
[x_{j,k-\widetilde{\nu}+1},x_{j,k+1}].
\)
In the even-odd partitioning, where $x_{j,k} = x_{j+1,2k}$, this means that
$2(k-\widetilde{\nu}+1) \leq \ell-\widetilde{\nu}+1$ and
$\ell+1\leq 2(k+1)$, leading to
$2k-\widetilde{\nu}+1 \leq \ell \leq 2k+1$.

The nonzero elements of $\mathbf{H}_j$ can be filled in by a loop over the
columns $k$, imposing two conditions.
The first condition is the normalisation
\(
\Phi_j(x) \vec{1}_j  = \omega_0(x) = \Phi_{j+1}(x) \vec{1}_{j+1},
\)
where we adopt the notation $\vec{1}_j$ for the vector of length $n_j$ whose
elements are all equal to one. Filling in the two-scale equation
(\ref{eq:2scale}) leads to $\Phi_{j+1}(x) \mathbf{H}_j \vec{1}_j
= \Phi_{j+1}(x) \vec{1}_{j+1}$, which can only hold for all possible $x$ if
\begin{equation}
\mathbf{H}_j \vec{1}_j = \vec{1}_{j+1}.
\label{eq:Hrowsum1}
\end{equation}
In a loop over the columns $k$, this expression can be used to find the
elements $H_{j;\ell,k}$ for the rows
$\ell \in \{2k-\widetilde{\nu}+1,2k-\widetilde{\nu}+2\} \cap \{0,1,\ldots\}$.
Indeed, for these values of $\ell$, we know that $H_{j;\ell,\kappa}=0$ if
$\kappa>k$, so
\begin{equation}
H_{j;\ell,k} = 1 - \sum_{\kappa=0}^{k-1} H_{j;\ell,\kappa}.
\label{eq:Hrowsum1bis}
\end{equation}

The second condition is the continuity of the $(\widetilde{\nu}-1)$th
derivative in the refinement knots.
Define the jumps of the $(\widetilde{\nu}-1)$th derivatives in the knots
$x_{j+1,i}$,
\[
\Delta_{j+1;i,\ell} = 
\lim_{x \to x_{j+1,i}^+}
{d^{\widehat{\nu}-1}\varphi_{j+1,\ell}(x)\over dx^{\widehat{\nu}-1}}
-
\lim_{x \to x_{j+1,i}^-}
{d^{\widehat{\nu}-1}\varphi_{j+1,\ell}(x)\over dx^{\widehat{\nu}-1}}.
\]
Then at the coarser scale, we need for $i\in O_{j+1}$ that
\[
\lim_{x \to x_{j+1,i}^+}
{d^{\widehat{\nu}-1}\varphi_{j,k}(x)\over dx^{\widehat{\nu}-1}}
-
\lim_{x \to x_{j+1,i}^-}
{d^{\widehat{\nu}-1}\varphi_{j,k}(x)\over dx^{\widehat{\nu}-1}}
=
0.
\]
Writing $\mathbf{\Delta}_{j+1;o}$ the
$(n_{j+1}-n_j)\times (n_{j+1}+\widetilde{\nu}-2)$
matrix with elements $\Delta_{j+1;i,\ell}$
for $i \in O_{j+1}$, this amounts to
\begin{equation}
\mathbf{\Delta}_{j+1;o}\mathbf{H}_j = \mathbf{0}.
\label{eq:Hjumps0}
\end{equation}
Not only the matrix $\mathbf{H}_j$ of unknowns has a band structure, also the
matrix $\mathbf{\Delta}_{j+1;o}$ has a limited number of nonzeros on each row.
Defining the index set
\(
\overline{S}_{j+1,\ell}
=
\{0,1,\ldots,n_{j+1}-1\} \cap \{\ell-\widetilde{\nu}+1,\ldots,\ell+1\}
=
S_{j+1,\ell} \cup \{\ell+1\},
\)
with $S_{j+1,\ell}$ as defined in (\ref{eq:defBB}), we have
$\Delta_{j+1;i,\ell} = 0$ for $i \not\in \overline{S}_{j+1,\ell}$. 
Hence, a row in $\mathbf{\Delta}_{j+1;o}$ involves a nonzero in
the $k$th column of $\mathbf{H}_j$ if $i \in O_{j+1}$ and if there exists an
index $\ell \in \{2k-\widetilde{\nu}+1,\ldots,2k+1\}$ so that
$\ell-\widetilde{\nu}+1 \leq i \leq \ell+1$.  This amounts to all
$i \in O_{j+1}$ with $2k-2\widetilde{\nu}+2\leq i \leq 2k+2$. Still working
with an even-odd split, we find $\widetilde{\nu}+1$ of such indices $i$,
each leading to an equation
\[
\sum_{\ell=2k-\widetilde{\nu}+1}^{2k+1} \Delta_{j+1;i,\ell} H_{j;\ell,k} = 0.
\]
The $\widetilde{\nu}+1$ equations involve $\widetilde{\nu}+1$ nonzeros in the
$k$-th column of $\mathbf{H}_j$, two of which have already been found by
(\ref{eq:Hrowsum1bis}).
The existence of a bandlimited refinement matrix guarantees the overcomplete
set to have an exact solution.

The boundary conditions (\ref{eq:BBleftboundarycond}) provide additional
conditions on the refinement matrix. Indeed, denoting by $\delta_{i,j}$ the
Kronecker delta (i.e., the $(i,j)$ element of the identity matrix), plugging
in $\varphi_{j+1,k}(x_0) = \delta_{0,k} = \varphi_{j,k}(x_0)$ into the
two-scale equation (\ref{eq:2scale}), leads to $H_{j;0,k} = \delta_{0,k}$.
Furthermore, by taking first and higher-order derivatives in (\ref{eq:2scale}),
one can also find that
\begin{equation}
H_{j;\ell,\ell+\kappa} = 0 \mbox{ for }\kappa = 1,2,\ldots.
\label{eq:Hboundary}
\end{equation}
As a result, the number of nonzeros in column $k$ of $\mathbf{H}_j$
equals $k+1$ for $k\leq \widehat{\nu}$. A similar reduced number of
nonzeros holds at the other boundary.

\section{Broken basis wavelets}
\label{sec:BBwavelets}

For the construction of a wavelet transform on top of a given refinement, the
refinement matrix can be factored into a sequence of bidiagonal lifting steps
\cite[Chapter~4]{jansen22:waveletbook}.
The factoring does not involve the boundary scaling functions. With $n_{j+1}$
and $n_j$ knots at resolution levels $j+1$ and $j$, the $n_{j+1}-n_j$ wavelet
functions can be constructed with $n_{j+1}$ fine and $n_j$ coarse scaling
functions. In particular, let $\mathbf{H}_{j;eo}$ be the $n_{j+1}\times n_j$
submatrix of $\mathbf{H}_j$, obtained by taking out the first
$\ceil{\widetilde{\nu}/2}-1$ and the last $\floor{\widetilde{\nu}/2}-1$ rows
and columns. Furthermore, let $\mathbf{H}_{j;e}$ be the ``even'' submatrix of
$\mathbf{H}_{j;eo}$, consisting of the even rows (and all columns) of
$\mathbf{H}_{j;eo}$. Similarly, $\mathbf{H}_{j;o}$ denotes the ``odd''
submatrix.

Then, a Euclid's algorithm \citep{daubechies98:factorlifting} can be proven to
factor the refinement according to one of two schemes, depending on the number
of nonzeros in each column of $\mathbf{H}_{j;eo}$, i.e., on the parameter
$\widetilde{\nu}$. The first scheme is
\begin{equation}
\left[\begin{array}{ll}\mathbf{H}_{j;e} & \mathbf{G}_{j;e}^{[0]} \\ 
                       \mathbf{H}_{j;o} & \mathbf{G}_{j;o}^{[0]}
\end{array}\right]
=
\left(\prod_{s=1}^{u}
\left[\begin{array}{ll}
\mathbf{I}_{n_j} & -\mathbf{U}_j^{[s]} \\ \mathbf{0} & \mathbf{I}_{n_j'}
\end{array}\right]
\left[\begin{array}{ll}
\mathbf{I}_{n_j} & \mathbf{0} \\ \mathbf{P}_j^{[s]} & \mathbf{I}_{n_j'}
\end{array}\right]\right)
\left[\begin{array}{ll}
\mathbf{D}_j & \mathbf{0} \\ \mathbf{0} & \mathbf{I}_{n_j'}
\end{array}\right].
\label{eq:factoring0}
\end{equation}
The alternative is given by
\begin{equation}
\left[\begin{array}{ll}\mathbf{H}_{j;e} & \mathbf{G}_{j;e}^{[0]} \\ 
                       \mathbf{H}_{j;o} & \mathbf{G}_{j;o}^{[0]}
\end{array}\right]
=
\left[\begin{array}{ll}
\mathbf{I}_{n_j} & \mathbf{0} \\ \mathbf{P}_j^{[0]} & \mathbf{I}_{n_j'}
\end{array}\right]
\left(\prod_{s=1}^{u}
\left[\begin{array}{ll}
\mathbf{I}_{n_j} & -\mathbf{U}_j^{[s]} \\ \mathbf{0} & \mathbf{I}_{n_j'}
\end{array}\right]
\left[\begin{array}{ll}
\mathbf{I}_{n_j} & \mathbf{0} \\ \mathbf{P}_j^{[s]} & \mathbf{I}_{n_j'}
\end{array}\right]
\right)
\left[\begin{array}{ll}
\mathbf{D}_j & \mathbf{0} \\ \mathbf{0} & \mathbf{I}_{n_j'}
\end{array}\right].
\label{eq:factoring1}
\end{equation}
In these equations, the number $u$ is given by
$u = \floor{(\widetilde{\nu}+1)/4}$.
The former factoring (\ref{eq:factoring0}) applies if the binary value
$r = \ceil{\widetilde{\nu}/2}-2u$ equals 0, while the latter
(\ref{eq:factoring1}) applies otherwise. 
The factorings adopt the notation $n_j' = n_{j+1}-n_j$ as index in some of the
identity matrices appearing in the expressions above. 
In general, identity matrices are denoted by $\mathbf{I}_m$, where $m$ stands
for the size of the matrix.
Both factorings consist of alternating matrix products, using two types of
matrices.
The first type involves $n_j'\times n_j$ bidiagonal matrices
$\mathbf{P}_j^{[s]}$, termed dual lifting or prediction matrices
(for reasons of interpretations that are beyond the scope of this paper).
The second type involves $n_j\times n_j'$ bidiagonal matrices
$\mathbf{U}_j^{[s]}$, termed primal lifting or update matrices.
The only difference between the two schemes (\ref{eq:factoring0}) and
(\ref{eq:factoring1}) is that the leftmost factor in the
first scheme is a prediction, while in the second scheme, the leftmost factor
is an update. The factoring also provides an $n_j\times n_j$ diagonal matrix
$\mathbf{D}_j$, which as a rescaling effect.

While (\ref{eq:factoring0}) and (\ref{eq:factoring1}) provide a factoring of
$\mathbf{H}_{j;eo}$ in the first place, it also defines the even and odd rows
of a matrix $\mathbf{G}_{j;eo}^{[0]}$, which can be used in the definition of
a primitive wavelet basis
\begin{equation}
\Psi_j^{[0]}(x) = \Phi_{j+1;eo}(x)\mathbf{G}_{j,eo}^{[0]}.
\end{equation}
In this expression, $\Phi_{j+1,eo}(x)$ are the interior scaling functions
within $\Phi_{j+1}(x)$, i.e., the scaling functions $\varphi_{j+1,k}(x)$ with
exception of the first $\ceil{\widetilde{\nu}/2}-1$ and the last
$\floor{\widetilde{\nu}/2}-1$ values of $k$.
Any other wavelet basis compatible with $\mathbf{H}_j$ can be found by one
single final update step \citep[Theorem~4.1.7]{jansen22:waveletbook}, i.e.,
\begin{equation}
\Psi_j(x) = \Phi_{j+1;eo}(x)\mathbf{G}_{j;eo},
\label{eq:wavelet}
\end{equation}
where the detail matrix $\mathbf{G}_{j;eo}$ follows from
\begin{equation}
\mathbf{G}_{j;eo} = \mathbf{G}_{j;eo}^{[0]} - \mathbf{H}_{j;eo}\mathbf{U}_j.
\label{eq:Ufinal}
\end{equation}
The $n_j\times n_j'$ matrix $\mathbf{U}_j$ can be used to equip the wavelet
basis $\Psi_j(x)$ with properties of interest such as vanishing moments or
control of variance propagation in statistical applications. The design of the
final matrix $\mathbf{U}_j$ is not limited to a bidiagonal or any other sparse
structure.

\begin{figure}[t!]
\hspace*{-20mm}
\includegraphics[width=1.2\textwidth]{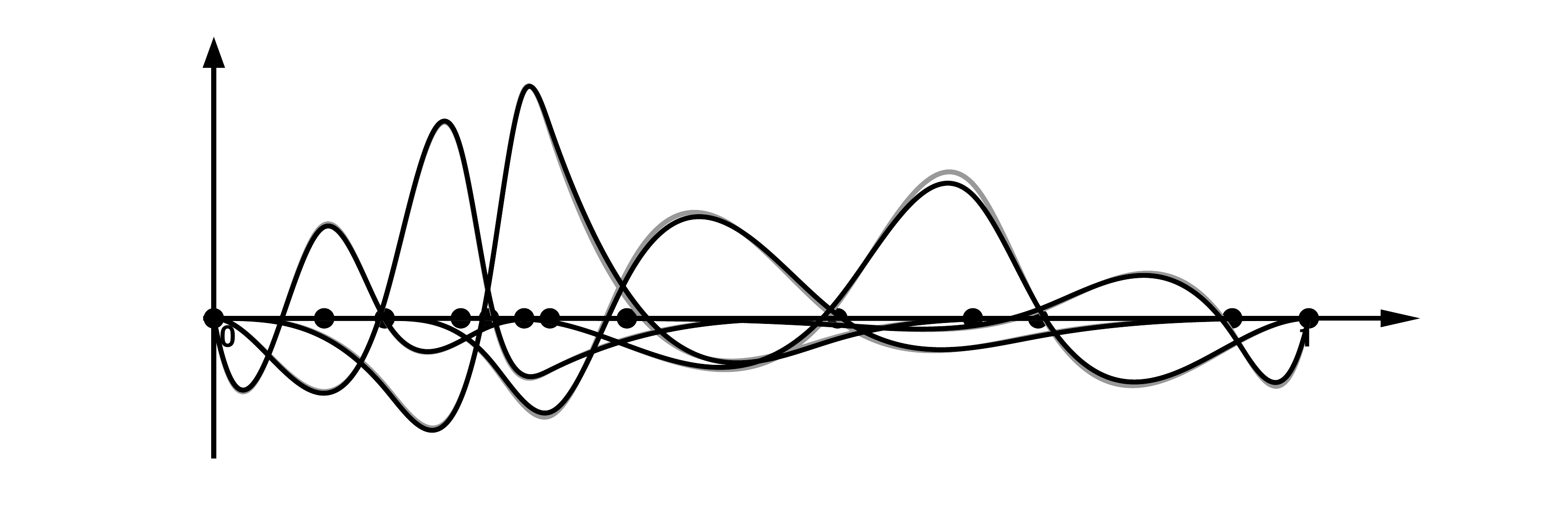}
\caption{Detail or wavelet functions for the refinement of the bases in Figure
\ref{fig:brokenbases}. The wavelet functions have two vanishing moments. Other
design options are possible and probably better in some applications. In grey:
cubic B-spline wavelets. In black (largely overlapping the grey curves),
wavelets composed of segments from
$\Omega(x) = \left[\begin{array}{cccc}1 & x & \sin(2\pi x) & \cos(2\pi x)
\end{array}\right]$.
}
\label{fig:BBwavelets}
\end{figure}

Figure \ref{fig:BBwavelets} presents the wavelets obtained for refinement of
the bases in Figure \ref{fig:brokenbases}. The wavelets have been designed to
have two primal vanishing moments. In other words, in this particular case, 
the final update matrix $\mathbf{U}_j$ has been taken to be bidiagonal, filling
in the two diagonals so that
\[
\int_0^1 \Psi_j(x) \,dx = \vec{0} = \int_0^1 \Psi_j(x) x\,dx.
\]

\section{The forward and inverse transforms}
\label{sec:BBWT}

\subsection{From bases to transforms}
\label{subsec:bases2WT}

The refinable and wavelet bases define a wavelet transform as follows.
Let $\vec{s}_{j+1}$ be fine scaling coefficients, then there exist coarse
scaling coefficients $\vec{s}_j$ and detail or wavelet coefficients
$\vec{d}_j$ so that
\[
\Phi_{j+1}(x)\vec{s}_{j+1} = \Phi_j(x)\vec{s}_j + \Psi_j(x)\vec{d}_j.
\]
Plugging in the two scale (\ref{eq:2scale}) and wavelet (\ref{eq:wavelet})
equations yields
\[
\Phi_{j+1}(x)\vec{s}_{j+1} = \Phi_{j+1}(x)\mathbf{H}_j\vec{s}_j +
\Phi_{j+1}(x)\mathbf{G}_j\vec{d}_j,
\]
where we have added $\widetilde{\nu}-2$ zero rows to $\mathbf{G}_{j;eo}$ for
the fine scaling functions on the two boundaries, so that
\(
\Phi_{j+1}(x)\mathbf{G}_j = \Phi_{j+1;eo}(x)\mathbf{G}_{j;eo}.
\)
Identification then leads to
\[
\vec{s}_{j+1} = \mathbf{H}_j\vec{s}_j + \mathbf{G}_j\vec{d}_j
=
\left[\begin{array}{cc} \mathbf{H}_j & \mathbf{G}_j\end{array}\right]
\left[\begin{array}{c} \vec{s}_j \\ \vec{d}_j\end{array}\right].
\]
Defining the transform matrices
$\widetilde{\mathbf{H}}_j$ and $\widetilde{\mathbf{G}}_j$ by
\begin{equation}
\left[\begin{array}{cc}
\widetilde{\mathbf{H}}_j & \widetilde{\mathbf{G}}_j\end{array}\right]^T
\left[\begin{array}{cc} \mathbf{H}_j & \mathbf{G}_j\end{array}\right]
=
\mathbf{I}_{n_{j+1}+\widetilde{\nu}-2}
=
\left[\begin{array}{cc} \mathbf{H}_j & \mathbf{G}_j\end{array}\right]
\left[\begin{array}{cc}
\widetilde{\mathbf{H}}_j & \widetilde{\mathbf{G}}_j\end{array}\right]^T,
\label{eq:defHtGt}
\end{equation}
we find one step of the multiresolution analysis or forward wavelet transform
as
\[
\left[\begin{array}{c} \vec{s}_j \\ \vec{d}_j\end{array}\right]
=
\left[\begin{array}{c} \widetilde{\mathbf{H}}_j^T \\
\widetilde{\mathbf{G}}_j^T\end{array}\right]
\vec{s}_{j+1}.
\]
The transform matrices $\widetilde{\mathbf{H}}_j$ and
$\widetilde{\mathbf{G}}_j$ follow from inverting
\(
\left[\begin{array}{cc} \mathbf{H}_j & \mathbf{G}_j\end{array}\right].
\)
In the lifting factoring, this inversion follows trivially, since
\[
\left[\begin{array}{ll}
\mathbf{I}_{n_j} & -\mathbf{U}_j^{[s]} \\ \mathbf{0} & \mathbf{I}_{n_j'}
\end{array}\right]^{-1}
=
\left[\begin{array}{ll}
\mathbf{I}_{n_j} & +\mathbf{U}_j^{[s]} \\ \mathbf{0} & \mathbf{I}_{n_j'}
\end{array}\right]
\mbox{ and }
\left[\begin{array}{ll}
\mathbf{I}_{n_j} & \mathbf{0} \\ \mathbf{P}_j^{[s]} & \mathbf{I}_{n_j'}
\end{array}\right]^{-1}
=
\left[\begin{array}{ll}
\mathbf{I}_{n_j} & \mathbf{0} \\ -\mathbf{P}_j^{[s]} & \mathbf{I}_{n_j'}
\end{array}\right].
\]

\subsection{Dealing with the boundaries}

The lifting factoring has been constructed on the interior basis functions.
Therefore, in order to incorporate the easy and fast inversion offered by this
factoring, we need to deal with the boundaries first. To this end, we write
\begin{equation}
\mathbf{H}_j = \left[\begin{array}{cc}
\mathbf{H}_{j;b,b} & \mathbf{H}_{j;b,eo} \\
\mathbf{H}_{j;eo,b} & \mathbf{H}_{j;eo}
\end{array}\right]
\mbox{ and }
\widetilde{\mathbf{H}}_j = \left[\begin{array}{cc}
\widetilde{\mathbf{H}}_{j;b,b} & \widetilde{\mathbf{H}}_{j;b,eo} \\
\widetilde{\mathbf{H}}_{j;eo,b} & \widetilde{\mathbf{H}}_{j;eo}
\end{array}\right],
\label{eq:defHHteob}
\end{equation}
where index $b$ refers to the $\widetilde{\nu}-2$ columns or rows on the
boundary. More precisely, $\mathbf{H}_{j;b,b}$ contains the elements
$H_{j;k,l}$ of $\mathbf{H}_j$ where
\(
k \in \{0,\ldots,\ceil{\widetilde{\nu}/2}-2\} \cap
      \{n_{j+1}+\ceil{\widetilde{\nu}/2}-2,\ldots n_{j+1}+\widetilde{\nu}-3\},
\)
and
\(
l \in \{0,\ldots,\ceil{\widetilde{\nu}/2}-2\} \cap
      \{n_j+\ceil{\widetilde{\nu}/2}-2,\ldots n_j+\widetilde{\nu}-3\}.
\)
The submatrix $\mathbf{H}_{j;eo}$ is the same as the one defined in Section
\ref{sec:BBwavelets}. It contains the interior rows and columns, to be
partitioned into even and odds for use in the factoring into lifting steps. The
submatrices $\mathbf{H}_{j;b,eo}$ and $\mathbf{H}_{j;eo,b}$ have boundary rows
on interior columns and interior rows on boundary columns respectively.
Expression (\ref{eq:defHHteob}) is a slight abuse of notation, in the sense
that it suggests that all boundary elements are on the upper left hand side of
the matrices $\mathbf{H}_j$ and $\widetilde{\mathbf{H}}_j$, whereas in reality,
the boundary elements are situated on both sides. The block matrices on the
right hand sides of the expressions thus tacitly include a rearrangement of
rows and columns.

Then the perfect reconstruction
of (\ref{eq:defHtGt}) implies that
\(
\widetilde{\mathbf{H}}_j^T\mathbf{H}_j = \mathbf{I}_{n_j+\widetilde{\nu}-2}.
\)
(Note that the rectangular matrices $\widetilde{\mathbf{H}}_j$ and
$\mathbf{H}_j$ do not satisfy $\mathbf{H}_j^T\widetilde{\mathbf{H}}_j =
\mathbf{I}_{n_{j+1}+\widetilde{\nu}-2}$.)
In terms of the submatrices defined in (\ref{eq:defHHteob}), this becomes
\begin{eqnarray}
\widetilde{\mathbf{H}}_{j;b,b}^T\mathbf{H}_{j;b,b}
+
\widetilde{\mathbf{H}}_{j;eo,b}^T\mathbf{H}_{j;eo,b}
& = &
\mathbf{I}_{\widetilde{\nu}-2},
\label{eq:HTHbb}
\\
\widetilde{\mathbf{H}}_{j;b,b}^T\mathbf{H}_{j;b,eo}
+
\widetilde{\mathbf{H}}_{j;eo,b}^T\mathbf{H}_{j;eo}
& = &
\mathbf{0}_{b,eo},
\label{eq:HTHbeo}
\\
\widetilde{\mathbf{H}}_{j;b,eo}^T\mathbf{H}_{j;b,b}
+
\widetilde{\mathbf{H}}_{j;eo}^T\mathbf{H}_{j;eo,b}
& = &
\mathbf{0}_{eo,b},
\label{eq:HTHeob}
\\
\widetilde{\mathbf{H}}_{j;b,eo}^T\mathbf{H}_{j;b,eo}
+
\widetilde{\mathbf{H}}_{j;eo}^T\mathbf{H}_{j;eo}
& = &
\mathbf{I}_{n_j}.
\label{eq:HTHeo}
\end{eqnarray}
From (\ref{eq:Hboundary}) it follows that $\mathbf{H}_{j;b,eo} =
\mathbf{0}_{b,eo}$. This means that the interior scaling functions can be
refined without intervention from the boundary functions, namely
\(
\Phi_{j;eo}(x) = \Phi_{j+1;eo} \mathbf{H}_{j;eo}.
\)
Substitution of $\mathbf{H}_{j;b,eo} = \mathbf{0}_{b,eo}$ into
(\ref{eq:HTHbeo}) leads to
\(
\widetilde{\mathbf{H}}_{j;eo,b}^T\mathbf{H}_{j;eo} = \mathbf{0}_{b,eo},
\)
which implies that
\(
\widetilde{\mathbf{H}}_{j;eo,b}^T = \mathbf{0},
\)
because the interior refinement $\mathbf{H}_{j;eo}$ is assumed to have full
rank $n_j$.
Further substitution of
\(
\widetilde{\mathbf{H}}_{j;eo,b}^T = \mathbf{0},
\)
into (\ref{eq:HTHbb}) reveals that
$\widetilde{\mathbf{H}}_{j;b,b}^T = \mathbf{H}_{j;b,b}^{-1}$.
In terms of forward wavelet transforms, this means that the coarse scaling
coefficients on the boundary $\vec{s}_{j,b}$ can be computed from a small
matrix inversion on the boundary,
\begin{equation}
\vec{s}_{j,b} = \mathbf{H}_{j;b,b}^{-1}\vec{s}_{j+1,b},
\label{eq:FWTb}
\end{equation}
not involving any interior scaling coefficient.
Substitution of the product of zero matrices
$\widetilde{\mathbf{H}}_{j;eo,b}^T\mathbf{H}_{j;b,eo}$ into
(\ref{eq:HTHeo}) also confirms that
\(
\widetilde{\mathbf{H}}_{j;eo}^T\mathbf{H}_{j;eo} = \mathbf{I}_{n_j}.
\)
The matrices $\widetilde{\mathbf{H}}_{j;eo}^T$ and $\mathbf{H}_{j;eo}$ follow
from the forward and inverse lifting schemes.
Unlike the squared boundary matrix
$\widetilde{\mathbf{H}}_{j;b,b}^T = \mathbf{H}_{j;b,b}^{-1}$,
the rectangular interior matrix $\widetilde{\mathbf{H}}_{j;eo}^T$ is not
uniquely defined for a given $\mathbf{H}_{j;eo}$.
It follows from the choice of the final update matrix used in th
design of the wavelet basis (\ref{eq:Ufinal}).

Finally, (\ref{eq:HTHeob}) can be written as
\[
\widetilde{\mathbf{H}}_{j;b,eo}^T
=
-\widetilde{\mathbf{H}}_{j;eo}^T\mathbf{H}_{j;eo,b}\mathbf{H}_{j;b,b}^{-1}.
\]
The interior coarse scaling coefficients are then given by
\[
\vec{s}_{j,eo}
=
\widetilde{\mathbf{H}}_{j;eo}^T \vec{s}_{j+1,eo}
+ \widetilde{\mathbf{H}}_{j;b,eo}^T \vec{s}_{j+1,b}
=
\widetilde{\mathbf{H}}_{j;eo}^T \left(\vec{s}_{j+1,eo}
- \mathbf{H}_{j;eo,b}\mathbf{H}_{j;b,b}^{-1}\vec{s}_{j+1,b}\right)
\]
Plugging in (\ref{eq:FWTb}), this becomes
\begin{equation}
\vec{s}_{j,eo}
=
\widetilde{\mathbf{H}}_{j;eo}^T \left(\vec{s}_{j+1,eo}-\mathbf{H}_{j;eo,b}
\vec{s}_{j,b}\right).
\label{eq:FWTseo}
\end{equation}
In other words, once the coarse scaling coefficients on the boundary,
$\vec{s}_{j,b}$, have been computed, the coarse scale interior coefficients
follow from a forward lifting scheme $\widetilde{\mathbf{H}}_{j;eo}^T$ applied
to slightly preprocessed fine scale interior coefficients
$\vec{s}_{j+1,eo}-\mathbf{H}_{j;eo,b}\vec{s}_{j,b}$.

The detail matrices can be written in two blocks
\[
\mathbf{G}_j = \left[\begin{array}{c} \mathbf{G}_{j;eo} \\ \mathbf{G}_{j;b}
\end{array}\right]
\mbox{ and }
\widetilde{\mathbf{G}}_j = \left[\begin{array}{c} 
\widetilde{\mathbf{G}}_{j;eo} \\ \widetilde{\mathbf{G}}_{j;b}
\end{array}\right],
\]
for which we know from Section \ref{subsec:bases2WT} that by construction
$\mathbf{G}_{j;b} = \mathbf{0}_{\widetilde{\nu}-2}$.
The perfect reconstruction (\ref{eq:defHtGt}) implies that
$\widetilde{\mathbf{G}}_j^T\mathbf{H}_j=\mathbf{0}$, which leads to
\[
\widetilde{\mathbf{G}}_{j;b} =
-\widetilde{\mathbf{G}}_{j;eo}\mathbf{H}_{j;eo,b}\mathbf{H}_{j;b,b}^{-1},
\]
and a similar result as in (\ref{eq:FWTseo}),
\begin{equation}
\vec{d}_j
=
\widetilde{\mathbf{G}}_{j;eo}^T \left(\vec{s}_{j+1,eo}-\mathbf{H}_{j;eo,b}
\vec{s}_{j,b}\right).
\label{eq:FWTdeo}
\end{equation}
In this expression the even and odd rows of $\widetilde{\mathbf{G}}_{j;eo}^T$
follow from the inverse of the factoring into lifting steps.

\subsection{Algorithmic overview}
\label{subsec:WTalgorithm}

All together, the implementation of a full forward wavelet transform in a
multiresolution analysis on nonequispaced knots requires the following steps
\begin{enumerate}
\item
Find the elements of the refinement matrix $\mathbf{H}_j$, following the steps
of Section \ref{sec:BBrefinement}.
\item
Factor the interior part of the refinement matrix $\mathbf{H}_{j;eo}$ into
lifting steps, defining a primitive detail matrix $\mathbf{G}_j^{[0]}$.
\item
Design an additional update $\mathbf{U}_j$ to define a wavelet basis through
the detail matrix $\mathbf{G}_j$ in (\ref{eq:Ufinal}). With $\mathbf{H}_j$ and
$\mathbf{G}_j$ chosen, the forward transform matrices
$\widetilde{\mathbf{H}}_j$ and $\widetilde{\mathbf{G}}_j$
are now fixed by (\ref{eq:defHtGt}). The following steps describe the
implementation of how to find these forward transform matrices.
\item
Find the course scaling coefficients on the boundaries, using (\ref{eq:FWTb}).
\item
Find the interior course scaling coefficients, using (\ref{eq:FWTseo}) and the
detail coefficients, using (\ref{eq:FWTdeo}).
Both expressions involve a matrix multiplication with
$\widetilde{\mathbf{H}}_{j;eo}^T$ and $\widetilde{\mathbf{G}}_{j;eo}^T$
respectively. These matrix multiplications are implemented by the running the
lifting operations from the factoring above.
\end{enumerate}

A few remarks conclude the design of the broken basis wavelet transform
\begin{enumerate}
\item
All steps in the scheme above are of linear complexity. This includes the
design of the refinement matrix $\mathbf{H}_j$ and its factoring into lifting
steps. As a result, the broken basis wavelet transform, operating on
nonequispaced knots, has the same order of complexity as the classical fast
wavelet transforms.
\item
The factoring in lifting steps is not unique. In particular, when applying the 
Euclid's algorithm to the refinement matrix of B-splines, the resulting lifting
scheme will slightly be different from the one constructed in
\citet{jansen16:noneqBsplinewavelets,jansen22:waveletbook}.
\end{enumerate}

\section{Discussion}
\label{sec:discussion}

The classical construction of wavelets concentrates on the two-scale equation
(\ref{eq:2scale}) on an equidistant set of knots. The two-scale equation
arrises within the framework of a multiresolution analysis, in which it is
imposed that all functions $\varphi_{j,k}(x)$ in the basis $\Phi_j(x)$ can be
written as normalised dilations and translations of a single father function,
i.e.,
\[
\varphi_{j,k} = 2^{j/2} \varphi(2^jx-k).
\]
On nonequispaced grids, the existence of a single father function is not
maintained, as illustrated by the case of B-spline scaling functions.
B-splines on equidistant knots, however, do have a single father function,
which can be found as the $\widetilde{\nu}$-fold convolution of the Haar
scaling basis (i.e., the indicator or characteristic function on $[0,1[$),
leading to the Cohen-Daubechies-Feauveau spline wavelets
\citep{coh-dau-fea92:biorthogonal}.

The existence of a single father function on equidistant knots is not
guaranteed by the refinable broken basis developed in the previous sections.
Indeed, the two-scale equation (\ref{eq:2scale}) does not impose the elements
of $\Phi_j(x)$ on an equidistant grid to be given by
\[
\varphi_{j,k}(x) = C\cdot\varphi_{j+1,k}(x/2)
\]
When the broken basis at finest scale is constructed by welding segments from
$\Omega(x)$, then the same set is used at all levels. As an example, if
$\Omega(x)$ includes the function $\sin(2\pi x)$, then the same sine function,
not a lower frequency version of it, appears in the segments of $\Phi_j(x)$ at
lower resolution levels $j<J$. In general, rescaled versions $\Omega(\alpha x)$
are not spanned by $\Phi_j(x)$ at any resolution level $j$, unless
$\Omega(\alpha x)$ happens to be a linear combination of $\Omega(x)$, as
is the case if $\Omega(x)$ contains the power functions.
Similar remarks hold for translations $\Omega(x+\beta)$ within one scale.

In may not be strictly necessary to have all dilations and translations of
$\Omega(x)$ within the space spanned by $\Phi_j(x)$. Indeed, the data to be
represented may show intermittent behaviour in scale or space, making a strict
condition of translation and dilation independence unrealistic. Nevertheless,
the example of taking $\omega_q(x) = \sin(2\pi x)$ as one the building
functions in $\Omega(x)$ illustrates the kind of problems that may arise from
not having any kind of dilation independence. Indeed, at coarse scales, the
basis $\Phi_j(x)$ is constructed from the same set $\Omega(x)$. With growing
intraknot distances, these segments start to capture the oscillating nature of
$\omega_q(x)$, leading to highly oscillating basis functions in $\Phi_j(x)$
at coarse scales.

As a conclusion, if we want the scaling functions at all resolution levels to
be dilations and translations of a single father function, then we need to
impose that
\begin{equation}
\Omega(\alpha x+\beta) = \Omega(x) A(\alpha,\beta),
\label{eq:diltrans}
\end{equation}
with $A(\alpha,\beta)$ independent from $x$.
The set of power functions, used in the construction of B-spline scaling bases,
satisfies these conditions. Taking the $k$-th derivative of (\ref{eq:diltrans})
w.r.t.~$x$ and letting $\alpha \to 0$, we obtain the condition
\begin{equation}
\Omega^{(k)}(\beta) = \Omega^{(k)}(x) {1 \over k!} {\partial^k \over
\partial\alpha^k} A(0,\beta).
\label{eq:diltransderiv}
\end{equation}
For (\ref{eq:diltransderiv}) to hold, we need that $\Omega^{(k)}(x)$ contains
the constant function $\omega_q(x)=1$ while the corresponding $q$-th row of
\(
{\partial^k \over \partial\alpha^k} A(0,\beta)
\)
equals $\Omega^{(k)}(\beta)$. Repeating the argument for $k=0,1,\ldots$ leads
to the conclusion that $\Omega(x)$ must contain the power functions
$\omega_q(x) = x^q$.
Further research is required to investigate the working of other sets of
$\Omega(x)$ not satisfying the strict conditions in (\ref{eq:diltrans}).

\section{Conclusion}
\label{sec:conclusion}

This paper has introduced the assembling of segments from a set of smooth
building functions into a refinable basis of scaling functions, followed by the
construction of a multiresolution wavelet transform. The resulting wavelets
have known smoothness and compact support, and they are defined on
nonequispaced knots on finite intervals, with special precautions in the
construction near the boundaries.
Unfortunetaly, taking arbitrary sets of building functions may result in
uncontrolled oscillations in the scaling basis, due to the lack of an imposed
father function generating all scaling functions by dilations and translations.
Ongoing and future research concentrates on controlling the lack of such a
father function.





\bibliographystyle{plainnat}
\bibliography{jansen23brokenbasispreprint}

\end{document}